\documentclass[11pt,a4paper]{article}
\usepackage[margin=27mm]{geometry}
\usepackage[T1]{fontenc}
\usepackage{lmodern,amsmath,amssymb,amsthm,mathtools,booktabs,microtype}
\usepackage[colorlinks=true,linkcolor=blue,citecolor=blue,urlcolor=blue]{hyperref}
\usepackage{bookmark}
\usepackage{cleveref}
\usepackage{authblk}
\usepackage{tikz}
\usetikzlibrary{arrows.meta,decorations.markings}
\newtheorem{theorem}{Theorem}
\newtheorem*{flowtheoremrestated}{Theorem~\ref{thm:flow}}
\newtheorem{lemma}[theorem]{Lemma}

\newtheorem{corollary}[theorem]{Corollary}
\theoremstyle{remark}
\newtheorem{remark}[theorem]{Remark}
\newcommand{\R}{\mathbb R}
\newcommand{\Z}{\mathbb Z}
\newcommand{\arc}[1]{\overrightarrow{#1}}
\newcommand{\larc}[1]{\overleftarrow{#1}}
\newcommand{\tauf}{\tau_f}
\newcommand{\cO}{\mathcal O}
\newcommand{\Ecal}{\mathcal E}

\hypersetup{pdftitle={Negative girth, circular coloring and cycle homomorphisms of planar graphs}}
\title{Jaeger conjecture and Jaeger-Zhang conjecture}

\author{Reza Naserasr}
\affil{Université Paris Cité, CNRS, F-75013 Paris, France. reza@irif.fr}

\date{}
\begin{document}
\maketitle

\begin{abstract}
We prove that every signed planar graph $\widehat G$ of negative girth at least
$\ell\geq 3$ satisfies $\chi_c(\widehat G)<\frac{2\ell}{\ell-2}$.
As a consequence, for every integer $k\geq 1$, every planar graph $G$ of
odd girth at least $4k+3$ satisfies $\chi_c(G)<\frac{2k+1}{k}$ and hence
admits a homomorphism to $C_{2k+1}$.

Our proofs build on recent results of Jafari on circular flows.
In particular, Jafari proved that every $(4k+2)$-edge-connected graph
has circular flow number strictly less than $2+\frac{1}{k}$.
We revisit his theorem, providing an alternative proof that emphasizes
the underlying graph-theoretic ideas. We also extend the result to
signed graphs and weaken the connectivity assumption by imposing
a lower bound only on the size of negative cuts.

The key step is a small-correction lemma associated with a distribution
on spanning trees. A splitting argument first gives the required
correction in an even-degree graph; doubling the edges removes the
degree-parity assumption. Plane duality then gives the coloring theorem.
We discuss tightness of bounds in the end.
\end{abstract}

\section{Introduction and Definitions}\label{sec:intro}

Homomorphisms from planar graphs to odd cycles have long been
motivated by Jaeger's conjecture and related questions about odd girth.
In signed graphs, negative cycles play the role of odd cycles, and
circular colorings quantify this perspective.\footnote{Although we work in the broader
setting of signed graphs, readers interested in ordinary graphs may assume
that every edge is negative, in which case negative cycles and negative cuts
are precisely odd cycles and odd cuts, respectively.} The dual notion is
circular flow. In tension language, coloring imposes constraints
along cycles, while flows enforce conservation at vertices.
Our main theorem establishes a strict upper bound for the circular
chromatic number of signed planar graphs in terms of their negative
girth. More precisely, if the negative girth is at least $\ell\ge3$,
then the circular chromatic number is strictly smaller than
$2\ell/(\ell-2)$.
The argument draws on recent work of Jafari on circular flows~\cite{Jafari},
which in turn uses the lattice-coset method of Alon, Buci\'c and
Davies~\cite{ABD26}.
We also return to his circular-flow theorem, giving an alternative
proof based on graph-theoretic ideas. At the heart of both
applications is a small-correction lemma derived from probability
distributions on spanning trees. In the planar setting, duality
allows us to translate the resulting flow statements into circular
coloring bounds.

Loops and parallel edges are allowed unless explicitly stated
otherwise. A loop contributes $2$
to the degree of its vertex and belongs to no cut. A cycle is a
connected graph all of whose vertices have degree $2$. A loop is a cycle of
length $1$, and two parallel edges form a cycle of length $2$.
A graph is \emph{even-degree} if every vertex has even degree,
without any requirement that it be connected. We use the term
\emph{Eulerian} for a connected even-degree graph.
A \emph{plane graph} is a planar graph with a fixed embedding in
the plane. For $X\subseteq V(G)$, the cut $\delta_G(X)$ consists
of the edges with exactly one endpoint in $X$, and we write
$d_G(X)=|\delta_G(X)|$. The notation $d_G(v)$ denotes the degree of
vertex $v$. Subscripts are omitted when the graph is clear,
and $G[X]$ denotes the subgraph induced by $X$.

A \emph{signed graph} $\widehat G=(G,\sigma)$ consists of a graph $G$
and a function $\sigma:E(G)\to\{+,-\}$. The sign of a walk is the
product of the signs of its edges, counted with multiplicity; its
length is its number of edges, also counted with multiplicity.
The \emph{negative girth} $g_-(\widehat G)$ is the length of a
shortest negative cycle, and is $\infty$ if no such cycle exists.
A signed graph without negative cycles is \emph{balanced}.
The \emph{odd girth} $g_o(G)$ is the length of a shortest odd cycle,
with the same convention for $\infty$.
Switching at a vertex reverses the signs of its incident nonloop
edges. More generally, switching at a set of vertices reverses the
signs of edges with exactly one endpoint in that set. Switching
preserves the sign of every closed walk. A \emph{homomorphism} between
signed graphs is a vertex map that, after switching the source if
necessary, maps each edge to an edge of the same sign in the target.
When no switching is allowed, we call the map a \emph{sign-preserving
homomorphism}.
For an ordinary graph $G$, the notation $(G,-)$ denotes the signed
graph with every edge negative; hence $g_o(G)=g_-((G,-))$.
For every integer $k\ge1$, the
notation $C_{-k}$ denotes a signed cycle of length $k$ whose sign
is negative. In particular, $C_{-1}$ is a negative loop and
$C_{-2}$ is a pair of parallel edges of opposite signs.

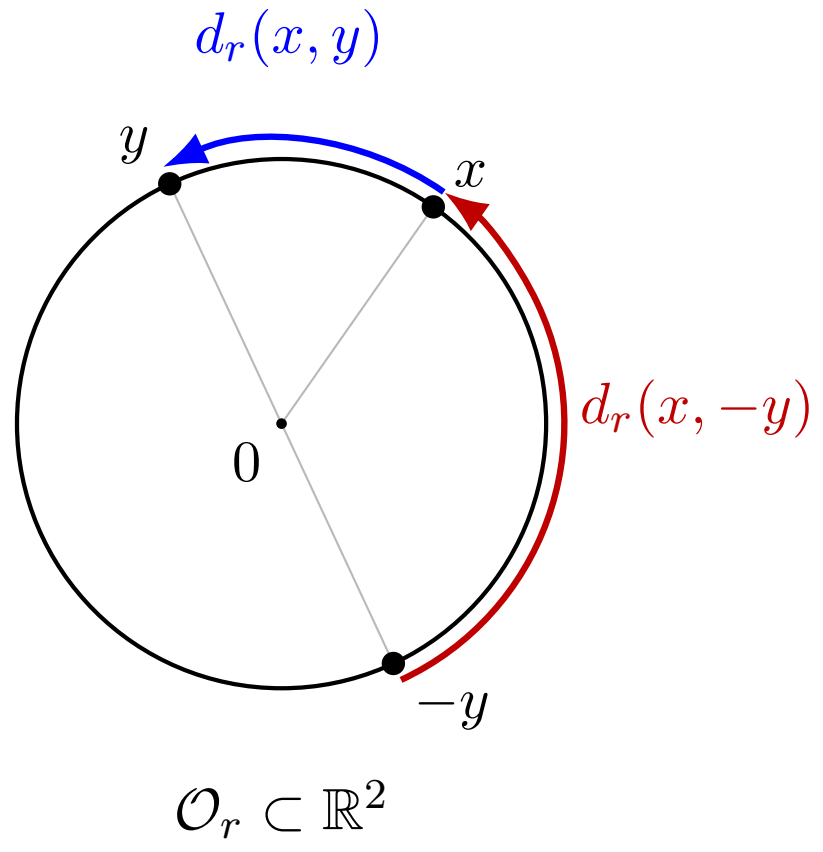
\begin{figure}[htbp]
\centering
\begin{tikzpicture}[scale=1.0,>=Latex]
 \def\rad{1.75}
 \def\outer{1.87}
 \coordinate (O) at (0,0);
 \coordinate (X) at (55:\rad);
 \coordinate (Y) at (115:\rad);
 \coordinate (NY) at (295:\rad);
 \draw[thick] (O) circle (\rad);
 \draw[gray!55] (O)--(X) (O)--(Y) (O)--(NY);
 \draw[blue,very thick,->] (55:\outer) arc[start angle=55,end angle=115,radius=\outer];
 \draw[red!75!black,very thick,->] (295:\outer) arc[start angle=295,end angle=415,radius=\outer];
 \fill (O) circle (1pt) node[below left] {$0$};
 \fill (X) circle (2.2pt) node[above right] {$x$};
 \fill (Y) circle (2.2pt) node[above left] {$y$};
 \fill (NY) circle (2.2pt) node[below right] {$-y$};
 \node[blue] at (0.05,2.55) {$d_r(x,y)$};
 \node[red!75!black] at (2.75,0.1) {$d_r(x,-y)$};
 \node at (0,-2.55) {$\mathcal O_r\subset\mathbb R^2$};
\end{tikzpicture}
\caption{Circular distances and antipodes.}
\label{fig:circular-distance}
\end{figure}

For $r\ge2$, let
\[
 \cO_r=\bigl\{x\in\R^2:\|x\|=\tfrac{r}{2\pi}\bigr\}
\]
denote the circle of circumference $r$, equipped with its
counterclockwise orientation. For $x,y\in\cO_r$, the distance
$d_r(x,y)$ is the length of a shortest arc joining them. The antipode
of $x$ is the point $-x$.
A \emph{circular $r$-coloring} of $\widehat G$ is a mapping
$\phi:V(G)\to\cO_r$ such that every edge $uv$ satisfies
\[
 d_r\bigl(\phi(u),-\sigma(uv)\phi(v)\bigr)\ge1.
\]
Thus, on a negative edge the two colors must be at distance at least
$1$, whereas on a positive edge one color must be at distance at
least $1$ from the antipode of the other; see \Cref{fig:circular-distance}.

The \emph{circular chromatic number} of a signed graph $\widehat G$,
denoted $\chi_c(\widehat G)$, is the
infimum of the feasible parameters $r\ge2$. Restricting to $(G,-)$
recovers the classical circular chromatic number of an ordinary
graph $G$. Positive loops impose no coloring restriction, whereas
a negative loop admits no circular coloring; thus
$\chi_c(C_{-1})=\infty$. This convention reverses the edge signs
in the original definition of Naserasr, Wang and Zhu~\cite{NWZ21}.
Switching at a vertex corresponds to replacing its color by its
antipode.

A circular $r$-coloring may also be described as a circular $r$-tension.
Its edge values satisfy prescribed bounds, and on every oriented cycle,
the sum of the values on forward edges minus the sum of the values on
backward edges is congruent to zero modulo $r$.
The dual notion is a \emph{circular $r$-flow} in the mono-directed sense
of~\cite{LNWZ24}. More precisely, after fixing a reference
orientation, a circular $r$-flow assigns a value modulo $r$ to each
arc, with $\partial f(v)\equiv0\pmod r$ at every vertex. With values
represented in $[0,r)$, each negative edge has value in $[1,r-1]$,
while each positive edge has value in $[0,r/2-1]\cup[r/2+1,r)$.
Reversing an arc negates its value modulo $r$.
A circular tension satisfies a cycle condition modulo $r$, and a
circular flow satisfies a vertex condition modulo $r$.
The \emph{circular flow number}
$\Phi_c(\widehat G)$ is the infimum of parameters $r\ge2$ for
which a circular $r$-flow exists. The boundary convention is given
in Section~\ref{sec:lengths}.

We first state the flow bound. Its proof appears in
Section~\ref{sec:flows}, after the small-correction argument.

\begin{theorem}\label{thm:flow}
For every integer $q\ge3$, if every negative cut of a signed
graph $\widehat G$ has size at least $q$, then
\[
 \Phi_c(\widehat G)<\frac{2q}{q-2}.
\]
\end{theorem}

Plane duality turns the flow bound into our main coloring theorem.

\begin{theorem}\label{thm:main}
For every integer $\ell\ge3$, every signed planar graph $\widehat G$
with $g_-(\widehat G)\ge\ell$ satisfies
\[
 \chi_c(\widehat G)<\frac{2\ell}{\ell-2}.
\]
\end{theorem}

The following consequences concern ordinary planar graphs, signed
bipartite planar graphs, and homomorphisms to negative cycles.

\begin{corollary}\label{cor:ordinary-color}
For every integer $k\ge1$, every planar graph $G$ with
$g_o(G)\ge2k+1$ satisfies
\[
 \chi_c(G)<\frac{2(2k+1)}{2k-1}
 =2+\frac{2}{k-\frac12}.
\]
\end{corollary}

\begin{corollary}\label{cor:bipartite-color}
For every integer $k\ge1$, every signed bipartite planar graph
$\widehat G$ with $g_-(\widehat G)\ge4k$ satisfies
\[
 \chi_c(\widehat G)<\frac{4k}{2k-1}.
\]
\end{corollary}

\begin{corollary}\label{cor:ordinary-target}
For every integer $k\ge1$, if $G$ is a planar graph of odd
girth at least $4k+3$, then
\[
 G\longrightarrow C_{2k+1}.
\]
\end{corollary}

\begin{corollary}\label{cor:bipartite-target}
For every integer $k\ge1$, if $\widehat G$ is a signed bipartite
planar graph of negative girth at least $4k$, then
\[
 \widehat G\longrightarrow C_{-2k}.
\]
\end{corollary}

The homomorphism corollaries follow from the circular-coloring
characterizations recalled in Section~\ref{sec:lengths}.

\section{Circular coloring and winding}\label{sec:lengths}

We consider a signed graph $\widehat G=(G,\sigma)$ and fix an arbitrary
orientation of $\widehat G$, denoted by $\overrightarrow{\widehat G}$.
An edge assignment is a function
\[
 \psi:E(\overrightarrow{\widehat G})\longrightarrow\mathbb R.
\]
If an edge is traversed against its chosen orientation, we take the
negative of its assigned value. Thus, for either traversal of an edge $uv$,
$\psi(\larc{uv})=-\psi(\arc{uv})$.
For any walk $W=(u_0,u_1,\ldots,u_\ell)$, whose vertices are not
necessarily distinct, we define
\[
 \psi(W)=\sum\limits_{i=0}^{\ell-1}\psi(\arc{u_i u_{i+1}}).
\]
When parallel edges occur, each step of $W$ specifies which edge is
traversed, and $\psi(\arc{u_i u_{i+1}})$ refers to that edge.
For an edge assignment $\psi$ and any
$X\subseteq V(G)$, its \emph{boundary} is the net value across the
cut $\delta_G(X)$:
\[
 \partial\psi(X)=
 \sum\limits_{\substack{uv\in E(G)\\u\in X,\ v\notin X}}
 \psi(\arc{uv}),
\]
where each edge of the cut is counted once in the direction from
$X$ to its complement. For a vertex $v$, we write
$\partial\psi(v)=\partial\psi(\{v\})$. Loops do not contribute
to the boundary of any vertex set. An assignment is a \emph{flow}
if $\partial\psi(v)=0$ for every vertex $v$. Throughout, flows on
signed graphs are understood in the mono-directed sense
of~\cite{LNWZ24}, rather than the usual bidirected sense.

Next, we express circular colorings in terms of arc lengths. Given a
circular $r$-coloring $\phi$, assign to each oriented edge $\arc{uv}$
the signed length $h(\arc{uv})$ of the shortest arc from $\phi(u)$ to
$\phi(v)$ on a positive edge, and to $-\phi(v)$ on a negative
edge, measured positively counterclockwise and negatively clockwise. The coloring condition gives
\[
 |h(\arc{uv})|\le \frac r2-1.
\]
In particular, this shortest arc is unique. As we traverse a cycle
$C$, the sum of these signed arc lengths, together with a half-turn
for each negative edge, is a multiple of $r$. Consequently,
\[
 h(C)\equiv\frac r2\,|E^-(C)|\pmod r.
\]
This identity may be viewed as a winding-number constraint;
see~\cite{NPW26} for a general treatment of winding numbers
in circular coloring.

The following edge-assignment formulation is the standard cycle-consistency
criterion for circular colorings, written in the sign convention used here.
We include its direct graph-theoretic proof.

\begin{lemma}\label{lem:short}
A signed graph has a circular $r$-coloring if and only if it has
an edge assignment $h$ satisfying
\begin{align}
 |h(e)|&\le \frac r2-1 &&(e\in E(G)),\label{eq:short-bound}\\
 h(C)&\equiv \frac r2\,|E^-(C)|\pmod r
 &&\text{for every cycle $C$, traversed in either direction.}
 \label{eq:cycle-parity}
\end{align}
Here $E^-(C)$ denotes the set of negative edges of $C$.
\end{lemma}
\begin{proof}
For a circular $r$-coloring, the shortest signed arc lengths
introduced above give an edge assignment. Its absolute value is at most
$r/2-1$. As we traverse a cycle, the signed arc lengths, together
with a half-turn for each negative edge, sum to a multiple of $r$.
This gives~\eqref{eq:cycle-parity}.

Conversely, suppose that $h$ satisfies the two conditions. For each
oriented edge $\arc{uv}$, consider the signed arc length
$h(\arc{uv})$ if $uv$ is positive, and $h(\arc{uv})+r/2$ if $uv$
is negative, taken modulo $r$. The signed arc length in the reverse
direction is its negative modulo $r$. By~\eqref{eq:cycle-parity},
these lengths sum to a multiple of $r$ around every cycle, traversed
in either direction. In each connected component, a starting vertex is assigned an
arbitrary point of $\cO_r$. The colors of the remaining vertices
are obtained by adding these signed arc lengths,
including the half-turns on negative edges, along a path from the
starting vertex, modulo $r$. The cycle condition guarantees that
the resulting colors are independent of the chosen paths.
The bound~\eqref{eq:short-bound} says that the colors
on a positive edge are at circular distance at most $r/2-1$,
while on a negative edge the color of one endpoint is at distance
at most $r/2-1$ from the antipode of the other. Both are precisely
the required circular-coloring inequalities.
\end{proof}

For every finite signed graph without negative loops, the infimum defining
its circular chromatic number is attained~\cite{NWZ21}.

The next two lemmas record known homomorphism characterizations.
The ordinary-graph characterization follows from classical circular
coloring theory (see~\cite{NWZ21}), while its signed bipartite analogue
is recorded in~\cite{NW23}. We state them in our sign convention.

\begin{lemma}\label{lem:ordinary-target}
For every graph $G$ and integer $k\ge1$,
\[
 G\longrightarrow C_{2k+1}
 \quad\Longleftrightarrow\quad
 \chi_c(G)\le\frac{2k+1}{k}.
\]
\end{lemma}

\begin{lemma}\label{lem:signed-target}
For a signed bipartite graph $\widehat G$ and $k\ge1$,
\[
 \widehat G\longrightarrow C_{-2k}
 \quad\Longleftrightarrow\quad
 \chi_c(\widehat G)\le\frac{4k}{2k-1}.
\]
\end{lemma}

\section{A small correction from spanning trees}\label{sec:correction}

A spanning tree is a connected acyclic subgraph containing every
vertex; a spanning forest is an acyclic subgraph containing every
vertex. A probability distribution on the spanning trees of a connected loopless
graph $G$ is a family of nonnegative numbers $(\lambda_T)$ satisfying
$\sum\limits_T\lambda_T=1$. The \emph{inclusion probability} of an edge
$e$ is $p(e)=\sum\limits_{T\ni e}\lambda_T$. The events that different
edges belong to the randomly chosen spanning tree need not be independent.
As in Section~\ref{sec:lengths}, we fix an arbitrary reference orientation
$\vec G$ and use the reversal convention for edge assignments.

For any assignment $\theta:E(G)\to\R$ and any set of edges
$A\subseteq E(G)$, we use the additive and multiplicative extensions
\[
 \theta^{+}(A)=\sum\limits_{e\in A}\theta(e),\qquad
 \theta^{\times}(A)=\prod\limits_{e\in A}\theta(e).
\]
In particular, these conventions apply when $A$ is the edge set of a
spanning tree or forest.

\begin{theorem}\label{thm:correction}
Suppose $G$ is a connected loopless graph equipped with a probability
distribution $(\lambda_T)$ on its spanning trees, and $p(e)$ denotes the
corresponding edge inclusion probabilities. For every assignment
$\varphi:E(\vec G)\to\R$, there exists an assignment
$h:E(\vec G)\to\R$ such that
\[
 |h(e)|\le\frac{p(e)}{2}\quad(e\in E(G)),\qquad
 \partial(\varphi+h)(v)\in\Z\quad(v\in V(G)).
\]
\end{theorem}

The correction principle underlying this theorem appears in the work
of Jafari~\cite{Jafari}, which builds on the
lattice-coset method of Alon, Buci\'c and Davies~\cite{ABD26}. The present formulation uses a
probability distribution on spanning trees. The proof of
Theorem~\ref{thm:correction} uses only spanning trees and edge
assignments. We organize the argument around two auxiliary lemmas.
The resulting prescribed-boundary corollary is the form used later
in the paper.

\begin{lemma}\label{lem:product}
Suppose $G$ is a connected loopless graph equipped with a probability
distribution $(\lambda_T)$ on its spanning trees, and let
$p(e)=\sum\limits_{T\ni e}\lambda_T$ denote the corresponding inclusion
probabilities. For every $\eta>0$, there exist positive edge weights
$\gamma(e)$ (depending on $\eta$) such that the distribution
\[
 Z=\sum\limits_T\gamma^{\times}(T),\qquad
 \lambda^\gamma_T=\frac{\gamma^{\times}(T)}Z
\]
has inclusion probabilities $p_\gamma(e)\le p(e)+\eta$.
\end{lemma}
\begin{proof}
We consider an assignment $\theta:E(G)\to\R$ satisfying
$\theta(e)\ge0$ for every edge $e$. The function $F$ is
defined on the set of all such assignments by
\[
 F(\theta)=\log\!\left(\sum\limits_T
       \frac{1}{\exp(\theta^{+}(T))}\right)
       +\sum\limits_e(p(e)+\eta)\theta(e).
\]
We first show that $F$ attains a minimum on its domain. By the
definition of the inclusion probabilities,
\[
 \sum\limits_T\lambda_T\theta^{+}(T)
 =\sum\limits_e p(e)\theta(e).
\]
Since $\sum\limits_T\lambda_T=1$, there is a spanning tree $T_0$ such that
\[
 \theta^{+}(T_0)\le\sum\limits_e p(e)\theta(e).
\]
All the terms in the sum inside the logarithm are positive. Retaining
only the term corresponding to $T_0$, we obtain
\[
\begin{aligned}
 F(\theta)
 &\ge -\theta^{+}(T_0)+\sum\limits_e(p(e)+\eta)\theta(e)\\
 &\ge \eta\sum\limits_e\theta(e).
\end{aligned}
\]
Consequently, the sublevel set of assignments satisfying
$F(\theta)\le F(0)$ is closed and bounded. Since $F$ is continuous,
it attains a minimum on this compact set. Because the zero assignment
belongs to the set and every assignment outside it has value greater
than $F(0)$, this is also a minimum on the whole domain.

A minimizing assignment is denoted by $\theta^*$, with
$\gamma(e)=\exp(-\theta^*(e))$. Differentiating $F$ with respect to
the coordinate $\theta(e)$, we obtain
\[
 D_eF(\theta^*)=p(e)+\eta-p_\gamma(e).
\]
If $\theta^*(e)>0$, this derivative is zero, since the coordinate can
be varied in either direction. If $\theta^*(e)=0$, its right derivative
is nonnegative, since increasing the coordinate is an allowed
variation at a minimum. Thus $D_eF(\theta^*)\ge0$ for every edge,
which gives $p_\gamma(e)\le p(e)+\eta$, as required.
\end{proof}

\pagebreak
We consider a connected loopless graph $G$ with a positive
edge-weight assignment $\gamma:E(G)\to\R_{>0}$ and
the corresponding probability distribution
$(\lambda_T^\gamma)$ on its spanning trees.
For distinct vertices $u,v$, the notation $P_T(u,v)$ denotes the unique
$u$--$v$ path in $T$, traversed from $u$ to $v$.
Its associated edge assignment
$\psi_T^{uv}:E(\vec G)\to\R$ takes value $1$ on edges of
$P_T(u,v)$ traversed in their reference direction, $-1$ on
edges traversed against that direction, and $0$ elsewhere.
The weighted average of these assignments is
\[
 J_{uv}=\sum\limits_T\lambda_T^\gamma\psi_T^{uv}.
\]
For an arbitrary edge assignment $\varphi:E(\vec G)\to\R$, define
its \emph{energy} with respect to $\gamma$ by
\[
 \Ecal_\gamma(\varphi)
 =\sum\limits_{e\in E(G)}\frac{\varphi(e)^2}{\gamma(e)}.
\]
If $\gamma(e)$ represents the conductance of an edge, then
$1/\gamma(e)$ is its resistance. Interpreting $\varphi(e)$ as
the current through that edge, $\Ecal_\gamma(\varphi)$ is
the total power dissipated by the network. This explains the
terminology.

For two edge assignments $\varphi,\rho:E(\vec G)\to\R$,
the \emph{weighted inner product} is
\[
 \langle\varphi,\rho\rangle_\gamma
 =\sum\limits_{e\in E(G)}
 \frac{\varphi(e)\rho(e)}{\gamma(e)}.
\]
Thus $\Ecal_\gamma(\varphi)=\langle\varphi,\varphi\rangle_\gamma$.
Two assignments are \emph{orthogonal with respect to $\gamma$}
if this inner product is zero. In electrical terminology, a
potential-driven current is orthogonal to every flow
(with boundary zero at every vertex), as follows from Ohm's law
and Kirchhoff's laws.

The next lemma supplies the electrical-network identity needed in our
proof. We give the argument directly in terms of weighted spanning trees;
no electrical-network theorem is invoked as a black box.

\begin{lemma}\label{lem:unit}
For every pair of distinct vertices $u,v\in V(G)$, the assignment
$J_{uv}$ satisfies
\[
 \partial J_{uv}(t)=
 \begin{cases}
  1,&t=u,\\
  -1,&t=v,\\
  0,&t\notin\{u,v\}.
 \end{cases}
\]
Moreover, $J_{uv}$ is orthogonal with respect to $\gamma$
to every flow $\rho:E(\vec G)\to\R$; that is,
\[
 \langle J_{uv},\rho\rangle_\gamma=0.
\]
Furthermore, if $u$ and $v$ are adjacent, with
$e=\arc{uv}$, then
\[
 J_{uv}(e)=p_\gamma(e),\qquad
 \Ecal_\gamma(J_{uv})=\frac{p_\gamma(e)}{\gamma(e)}.
\]
\end{lemma}
\begin{proof}
For any spanning tree or spanning forest $F$, its weight is
\[
 \gamma(F)=\prod\limits_{e\in F}\gamma(e).
\]
The total weight of all spanning trees is
\[
 Z=\sum\limits_{T\text{ spanning tree}}\gamma(T).
\]
The edge weights $\gamma$ induce a probability distribution on
spanning trees, given by $\lambda_T^\gamma=\gamma(T)/Z$.

The family $\mathcal F_{uv}$ consists of all spanning forests with
exactly two components, one containing $u$ and the other containing $v$.
For each $F\in\mathcal F_{uv}$, the vertex sets of these components
are denoted by $U_F$ and $V_F$, respectively. These forests give rise
to a function $q_{uv}$ on $V(G)$, defined by
\[
 q_{uv}(x)=\frac1Z
 \sum\limits_{\substack{F\in\mathcal F_{uv}\\x\in U_F}}\gamma(F).
\]
For every edge $xy$, oriented from $x$ to $y$ according to the fixed
reference orientation, we have
\begin{equation}\label{eq:forest-gradient}
 J_{uv}(xy)=\gamma(xy)\bigl(q_{uv}(x)-q_{uv}(y)\bigr).
\end{equation}
Indeed, a spanning tree $T$ contributes to $J_{uv}(xy)$ precisely
when its $u$--$v$ path contains $xy$. Removing this edge produces a
spanning forest with exactly two components, separating both $u$
from $v$ and $x$ from $y$. Conversely, adding $xy$ to such a forest
produces a spanning tree of weight $\gamma(xy)\gamma(F)$.
The path contribution is positive when $x\in U_F$ and negative
when $y\in U_F$. Summing the corresponding tree weights proves~\eqref{eq:forest-gradient}.

Every path assignment $\psi_T^{uv}$ has boundary $1$ at $u$, $-1$
at $v$, and $0$ elsewhere. Since the coefficients $\lambda_T^\gamma$
sum to $1$, their weighted average $J_{uv}$ has the same boundary.

For every flow $\rho$, equation~\eqref{eq:forest-gradient}
expresses the weighted inner product of $J_{uv}$ and $\rho$ as
\[
\begin{aligned}
 \langle J_{uv},\rho\rangle_\gamma
 &=\sum\limits_{xy\in E(G)}
   \frac{J_{uv}(xy)\rho(xy)}{\gamma(xy)}\\
 &=\sum\limits_{xy\in E(G)}
   \bigl(q_{uv}(x)-q_{uv}(y)\bigr)\rho(xy)\\
 &=\sum\limits_{z\in V(G)}q_{uv}(z)\,\partial\rho(z)\\
 &=0.
\end{aligned}
\]
The third equality collects the edge contributions corresponding to
each vertex (summation by parts). The last equality holds because
$\rho$ is a flow, so its boundary is exactly zero at every vertex.
This establishes weighted orthogonality.

Finally, suppose $e=\arc{uv}$. The $u$--$v$ path in a spanning tree
$T$ contains $e$ if and only if $e\in T$, and in that case it
traverses $e$ from $u$ to $v$. Consequently,
\[
 J_{uv}(e)=\sum\limits_{\substack{T\text{ spanning tree}\\e\in T}}
 \lambda_T^\gamma=p_\gamma(e).
\]
Reducing the value of $J_{uv}$ by $1$ on $e$, while leaving its
values on all other edges unchanged, produces a flow, because the
boundary becomes zero at every vertex. Orthogonality of $J_{uv}$
to this flow gives
\[
 0=\Ecal_\gamma(J_{uv})-\frac{J_{uv}(e)}{\gamma(e)}.
\]
Therefore,
\[
 \Ecal_\gamma(J_{uv})
 =\frac{J_{uv}(e)}{\gamma(e)}
 =\frac{p_\gamma(e)}{\gamma(e)},
\]
as claimed.
\end{proof}

\begin{proof}[Proof of Theorem~\ref{thm:correction}]
We fix positive weights $\gamma$ and minimize $\Ecal_\gamma(h)$
subject to $\partial(\varphi+h)$ being integer-valued.
A minimum exists: the feasible set is nonempty (take $h=-\varphi$)
and closed, and the energy has compact sublevel sets.

For every flow $\rho$, the assignments $h+t\rho$ remain feasible
for all real $t$, because their boundaries do not change.
Since $h$ minimizes the energy, differentiating
$\Ecal_\gamma(h+t\rho)$ at $t=0$ gives
$\langle h,\rho\rangle_\gamma=0$.

For an edge $e=\arc{uv}$, Lemma~\ref{lem:unit} shows that
the assignment $J_{uv}$ has boundary $1$ at $u$, $-1$ at $v$,
and $0$ elsewhere. Consequently, the assignment obtained from
$J_{uv}$ by decreasing its value on $e$ by $1$, while leaving
all other values unchanged, is a flow. Since $h$ is orthogonal
to every flow, we obtain
\[
 \langle h,J_{uv}\rangle_\gamma-\frac{h(e)}{\gamma(e)}=0,
\]
and hence
\[
 \langle h,J_{uv}\rangle_\gamma=\frac{h(e)}{\gamma(e)}.
\]

The assignments $h+J_{uv}$ and $h-J_{uv}$ are also feasible:
the boundary of $\varphi+h$ is already integer-valued, and
adding or subtracting $J_{uv}$ changes it only by integers.
Since $h$ minimizes the energy, neither change can decrease it.
Using the energy identity in Lemma~\ref{lem:unit}, we obtain
\[
 0\le\Ecal_\gamma(h\pm J_{uv})-\Ecal_\gamma(h)
   =\frac{\pm2h(e)+p_\gamma(e)}{\gamma(e)}.
\]
The two choices of sign imply
$|h(e)|\le p_\gamma(e)/2$ for every edge $e$.

For each positive integer $j$, Lemma~\ref{lem:product}, with
$\eta=1/j$, provides positive edge weights $\gamma_j$ such that
\[
 p_{\gamma_j}(e)\le p(e)+\frac1j
\]
for every edge $e$. The argument above, applied with these weights,
gives an assignment $h_j$ for which
$\partial(\varphi+h_j)$ is integer-valued and
\[
 |h_j(e)|\le\frac{p(e)+1/j}{2}.
\]
In particular, the values $h_j(e)$ are bounded independently
of $j$. Since $G$ has finitely many edges, the sequence of
assignments $(h_j)$ has a convergent subsequence. We denote
its limit by $h$.

For each edge $e$, passing to the limit in the preceding
inequality gives
\[
 |h(e)|\le\frac{p(e)}2,
\]
as required.

It remains to verify that the boundary of $\varphi+h$ is
integer-valued. For every vertex $v$, the boundary
$\partial(\varphi+h_j)(v)$ is an integer. Moreover, along our
convergent subsequence, these boundaries converge to
$\partial(\varphi+h)(v)$, since the boundary at a vertex is a
finite sum of edge values. A convergent sequence of integers
has an integer limit. Consequently,
\[
 \partial(\varphi+h)(v)\in\Z
\]
for every vertex $v$. Thus $h$ satisfies both required
conditions, completing the proof.
\end{proof}

\begin{corollary}\label{cor:prescribed-boundary}
Suppose $G$ is a connected loopless graph equipped with a probability
 distribution $(\lambda_T)$ on its spanning trees, and $p(e)$ denotes
 the probability that a randomly chosen spanning tree contains $e$.
For every mapping $\beta:V(G)\to\R$ satisfying
\[
 \sum\limits_{v\in V(G)}\beta(v)\in\Z,
\]
there exists an assignment $h:E(\vec G)\to\R$ such that
\[
 \partial h(v)\equiv\beta(v)\pmod1\quad(v\in V(G)),
 \qquad |h(e)|\le\frac{p(e)}2\quad(e\in E(G)).
\]
\end{corollary}
\begin{proof}
The sum of the values of $\beta$ is an integer. Subtracting this
integer from the value at one vertex gives a function $\beta'$ with
sum zero and the same values modulo one.

We now choose a spanning tree rooted at an arbitrary vertex. Starting
at its leaves and proceeding toward the root, the edges receive values
so that the boundary at each nonroot vertex is the prescribed value
of $\beta'$. Equivalently, an edge directed from a subtree toward its
parent receives the sum of $\beta'$ over that subtree, with zero
assigned to edges outside the tree. Since the sum of $\beta'$ is zero,
the boundary at the root also has the required value. In this way we
obtain an assignment $w$ with $\partial w=\beta'$.

Theorem~\ref{thm:correction}, applied to $\varphi=-w$, gives an
assignment $h$ with $|h(e)|\le p(e)/2$ and
$\partial(h-w)(v)\in\Z$ at every vertex. Consequently
$\partial h(v)\equiv\partial w(v)\equiv\beta(v)\pmod1$,
as required.
\end{proof}

\section{Prescribing half-integral boundaries}\label{sec:parity}

Throughout this section, signed graphs are connected and have no
negative loops, unless stated otherwise. The \emph{sign of a cut} is the product of the signs of its
edges. A cut is \emph{positive} or \emph{negative} according as this
product is $+$ or $-$. We say that a cycle is \emph{inversed} when
the signs of all its edges are reversed. Inversing a cycle leaves the
sign of every cut unchanged, since every cut meets the cycle in an
even number of edges. For further discussion of inversing in signed
graphs, see Li, Naserasr, Wang and Zhu~\cite{LNWZ24}.

For a signed graph $\widehat G=(G,\sigma)$, let $d^-(v)$ denote the
number of negative edges incident with $v$, counted with multiplicity. We seek an assignment whose boundary satisfies
\begin{equation}\label{eq:T-boundary}
 \partial h(v)\equiv\frac{d^-(v)}2\pmod1
 \qquad(v\in V(G)).
\end{equation}
By the handshaking lemma, $\sum\limits_{v\in V(G)}d^-(v)$ is even.
Moreover, a cut $\delta_G(X)$ is negative if and only if
$\sum\limits_{v\in X}d^-(v)$ is odd.

We use two classical graph results.

\subsection{Fractional spanning-tree packing}\label{subsec:fractional-packing}

A \emph{fractional spanning-tree packing} assigns weights $w_R\ge0$ to the spanning
trees $R$, with $\sum\limits_{R\ni e}w_R\le1$ on each edge. Its maximum
total weight is denoted by $\tauf(F)$. The fractional form of the
Nash--Williams--Tutte theorem~\cite{NashWilliams61,Tutte61} gives,
for every connected loopless graph $F$ with at least two vertices,
\begin{equation}\label{eq:fractional-packing}
 \tauf(F)=\min_{\mathcal P:\,|\mathcal P|\ge2}
       \frac{|E_F(\mathcal P)|}{|\mathcal P|-1},
\end{equation}
where $\mathcal P$ ranges over partitions of $V(F)$ into nonempty
parts and $E_F(\mathcal P)$ consists of the edges between parts.
Normalizing a maximum packing gives a tree distribution with
$p(e)\le1/\tauf(F)$. Thus Corollary~\ref{cor:prescribed-boundary}
gives an assignment with any required boundary values modulo $1$
whose sum is an integer, and with
\begin{equation}\label{eq:packing-boundary}
 |h(e)|\le\frac1{2\tauf(F)}.
\end{equation}

\subsection{Splitting in even-degree signed graphs}\label{subsec:sign-splitting}

For a signed graph $\widehat G=(G,\sigma)$, extend the edge-sign
function to the vertices by setting
\[
\sigma(v)=\prod\limits_{e\in\delta_G(v)}\sigma(e)=(-1)^{d^-(v)},
\]
where a loop is counted twice in $d^-(v)$. Thus a vertex is positive
or negative according as $\sigma(v)$ is $+$ or $-$. For
$X\subseteq V(G)$, put $\sigma(X)=\prod\limits_{v\in X}\sigma(v)$.
Cancellation of the internal edge signs gives
\begin{equation}\label{eq:vertex-cut-sign}
	\sigma(\delta_G(X))=\sigma(X).
\end{equation}
In particular, inversing a cycle preserves every vertex sign and
every cut sign.

Suppose that $zx$ and $zy$ are distinct nonloop edges with the same
sign $s$. A \emph{splitting at $z$} replaces these two edges by an
edge $xy$ of sign $s$. This operation preserves every vertex sign,
and hence, by~\eqref{eq:vertex-cut-sign}, the sign of the cut
associated with every vertex set. If $x=y$, the resulting loop is
deleted. More generally, deleting a loop of either sign changes
neither the vertex signs nor any cut. We allow inversing before a
splitting.

\begin{lemma}\label{lem:signed-splitting}
	Let $q\ge 2$ be even, and let $\widehat G=(G,\sigma)$ be a connected,
	loopless even-degree signed graph in which every negative cut has size at
	least $q$. Let $z$ be a vertex such that either
	\begin{enumerate}
		\item $z$ is positive and $d_G(z)>0$, or
		\item $z$ is negative and $d_G(z)>q$.
	\end{enumerate}
	After inversing a cycle if necessary, a splitting can be performed
	at $z$ so that every negative cut of the resulting graph still has
	size at least $q$. All vertex signs are preserved, and the resulting
	graph is again even-degree.
\end{lemma}

\begin{proof}
	Consider two distinct nonloop edges $zx$ and $zy$ incident with $z$,
	possibly with $x=y$, without imposing any condition on their signs.
	Let $G'$ be the unsigned graph obtained from $G$ by deleting these
	two edges and adding an edge $xy$. Call the pair of edges
	\emph{admissible} if $d_{G'}(X)\ge q$ whenever $\sigma(X)=-$,
	where the sign is taken in the original signed graph $\widehat G$.
	For any $X\subseteq V(G)$, if $X$ separates $z$ from both $x$ and
	$y$, then $d_{G'}(X)=d_G(X)-2$; otherwise,
	$d_{G'}(X)=d_G(X)$. Since every cut of $G$ has even size,
	admissibility fails precisely when some negative cut of size $q$
	separates $z$ from both $x$ and $y$.
	
	A set $X\subseteq V(G)$ is called \emph{tight} if
	$\sigma(X)=-$ and $d_G(X)=q$. Thus the pair of edges $zx,zy$ is not
	admissible if and only if some tight set contains both $x$ and $y$
	but not $z$.
	
	\emph{Claim 1.} Let $zy$ be an edge. If $A$ and $B$ are two tight sets
	with $y\in A\cap B$ and $z\notin A\cup B$, then $A\cup B$ is tight.
	
	\emph{Proof of Claim 1.} Suppose first that $\sigma(A\cap B)=+$.
	Then $A\setminus B$ and $B\setminus A$ both have negative sign.
	Since $zy$ joins $A\cap B$ to $V(G)\setminus(A\cup B)$, we obtain
	\[
	2q=d_G(A)+d_G(B)
	=d_G(A\setminus B)+d_G(B\setminus A)
	+2e_G(A\cap B,V(G)\setminus(A\cup B))\ge 2q+2,
	\]
	a contradiction. Here $e_G(P,Q)$ counts the edges between disjoint
	sets $P,Q$. Thus $\sigma(A\cap B)=-$, which also implies
	$\sigma(A\cup B)=-$. By submodularity,
	\[
	2q=d_G(A)+d_G(B)
	\ge d_G(A\cap B)+d_G(A\cup B)\ge 2q.
	\]
	It follows that $d_G(A\cup B)=q$, and so $A\cup B$ is tight.
	\hfill$\diamond$
	
	\emph{Claim 2.} No tight set $X$ with $z\notin X$ contains all
	neighbors of $z$.
	
	\emph{Proof of Claim 2.} Suppose that a tight set $X$, not containing
	$z$, contains all neighbors of $z$. If $z$ is positive, then
	$\sigma(X\cup\{z\})=-$ and
	$d_G(X\cup\{z\})=q-d_G(z)<q$, a contradiction. If $z$ is
	negative, then $q=d_G(X)\ge d_G(z)>q$, again a contradiction.
	\hfill$\diamond$
	
	Choose an edge $zx$. If no tight set contains $x$ but not $z$, then
	pairing $zx$ with any other edge incident with $z$ gives an admissible
	pair. Otherwise, let $X$ be an inclusion-maximal tight set containing
	$x$ but not $z$. By Claim 2, there is a neighbor $y$ of $z$ outside
	$X$. If a tight set containing both $x$ and $y$ but not $z$ existed,
	then by Claim 1 its union with $X$ would be a larger tight set not
	containing $z$, contradicting maximality. Thus $zx,zy$ is admissible.
	
	It remains to arrange equal signs. Suppose that the chosen edges
	$e=zx$ and $f=zy$ have opposite signs. Since every vertex of $G$
	has even degree, $G$ has no bridges. If $G-\{e,f\}$ were
	disconnected, some cut of $G$ would consist exactly of $e,f$.
	This cut would be negative and would become empty after splitting,
	contrary to admissibility. Hence $G-\{e,f\}$ is connected.
	There is therefore a cycle containing $e$ and avoiding $f$.
	Inversing this cycle makes the signs of $e$ and $f$ equal, while
	preserving all vertex signs and the admissibility of the pair.
	
	Now split the pair, assigning its common sign to the new edge.
	Vertex signs are preserved, so admissibility gives the required
	bound on negative cuts. Degrees remain even, also after deletion
	of a resulting loop.
\end{proof}

\begin{corollary}\label{cor:signed-regular-reduction}
	Let $q\ge2$ be even, and let $\widehat G$ be an even-degree signed
	graph in which every negative cut has size at least $q$. By a
	sequence of inversings, splittings and loop deletions, one can
	obtain a signed graph on the same vertex set in which every
	originally positive vertex is isolated and every originally
	negative vertex has degree exactly $q$. Every vertex retains its
	sign, and every negative cut retains size at least $q$ throughout
	the sequence.
\end{corollary}

\begin{proof}
	Delete loops and apply Lemma~\ref{lem:signed-splitting}
	within the connected component containing a positive vertex of
	positive degree or a negative vertex of degree greater than $q$,
	deleting any new loop. The hypothesis on negative cuts passes to
	each connected component, and remains valid when components split,
	since every negative cut within a component is also a negative
	cut of the whole graph. Each step reduces
	the number of edges, so the process terminates. A negative vertex
	always has degree at least $q$, because its singleton cut is
	negative. The stated degrees therefore hold at termination.
\end{proof}

The classical antecedents of this reduction are the splitting-off
results of Lov\'asz~\cite{Lovasz76} for Eulerian graphs and
Mader~\cite{Mader78} for local edge-connectivity. Zhang's
splitting theorem~\cite{Zhang02} preserves odd-edge-connectivity
and does not require planarity. Szigeti~\cite[Theorem~8]{Szigeti08}
subsequently generalized Zhang's theorem. Its cut argument is
used above. Our proof also supplies the inversing step needed for
the signed formulation.

The plane counterparts are folding lemmas. Klostermeyer and
Zhang~\cite{KlostermeyerZhang00} proved a folding lemma preserving
odd girth. Naserasr, Rollov\'a and Sopena~\cite[Lemma~3.1]{NRS13}
established the signed bipartite version preserving negative
girth, which is precisely the dual of the present statement when $G$ is planar.

\subsection{General boundary theorem}\label{subsec:general-boundary}

We now use the splitting lemma to establish the boundary theorem
for general signed graphs.

\begin{theorem}\label{thm:general-boundary}
Let $\widehat G$ be a connected signed graph in which every negative
cut has size at least $q\ge2$. There is an edge assignment
$h$ satisfying
\[
\partial h(v)\equiv\frac{d^-(v)}2\pmod1
\quad(v\in V(G)) \quad \text{ and } \quad
|h(e)|<\frac1q\qquad(e\in E(G)).
\]
\end{theorem}
\begin{proof}
First, when $G$ is an even-degree graph and $q$ is even, choose a
counterexample minimizing $|V(G)|+|E(G)|$. Positive loops can be assigned
zero and removed, since they belong to no cut and do not affect the
boundary congruences. If there are no negative cuts, then $d^-(v)$ is even at every vertex,
and the zero assignment suffices. We may thus take $G$ loopless
and assume a negative cut exists.

If some vertex $z$ is positive, Lemma~\ref{lem:signed-splitting}
provides an admissible splitting at $z$, after inversing a cycle
if necessary. Inversing does not change any vertex sign and
therefore leaves the required boundary values
unchanged modulo $1$. Let $\widehat G'$ be the signed
graph after splitting and deleting any resulting loop. It has
even degrees, fewer edges, and the same vertex signs; every
negative cut still has size at least $q$. Each nontrivial
component has an assignment by minimality.
If the splitting replaces $zx,zy$ by an edge $xy$, its assigned
value can be carried along the path $xzy$, with orientations
chosen so that the boundary at $z$ is zero. The other boundary
values are unchanged. If $x=y$, the resulting loop is deleted
and we assign zero to the two original edges. In either case,
the prescribed boundary congruences and the strict bound on
edge values hold on the original graph, a contradiction.

Consequently, every vertex \(v\in V(G)\) is negative. Hence \(\delta_G(v)\) is negative, and therefore
\[
d_G(v)=|\delta_G(v)|\ge q.
\]

Writing $n=|V(G)|\ge2$, we have
\[
 |E(G)|\ge\frac{qn}{2}>\frac q2(n-1).
\]
Let $S$ be an inclusion-minimal nonempty vertex subset satisfying
$|E(G[S])|>\frac q2(|S|-1)$, and denote its induced subgraph by
$F=G[S]$. Then $|S|\ge2$, and for every partition
$\mathcal P$ of $S$ into $p\ge2$ nonempty parts, minimality implies
\[
 |E_F(\mathcal P)|
 =|E(F)|-\sum\limits_{X\in\mathcal P}|E(G[X])|
 >\frac q2(p-1).
\]
In particular $F$ is connected, and by the definition of $\tauf$ we have
$\tauf(F)>q/2$.

Contract the connected subgraph $F=G[S]$ to a vertex $s$, deleting its internal edges and
retaining the signs of the remaining edges. The quotient has even
degrees, and every negative cut of the quotient corresponds to a
negative cut of $\widehat G$ of the same size. Since $|S|\ge2$, the quotient has fewer vertices than $G$ and
no more edges. By the choice of $\widehat G$ as a counterexample
minimizing $|V(G)|+|E(G)|$, the quotient admits an assignment $h'$
satisfying the required boundary congruences and $|h'(e)|<1/q$
for every edge $e$ of the quotient.
Retain $h'$ on the edges outside $F$. For each $v\in S$,
let $\eta(v)$ be the contribution of these edges to the boundary at $v$.
The required boundary value modulo $1$ for the assignment inside $F$
is therefore
\[
 \beta(v)=\frac{d^-(v)}{2}-\eta(v)\qquad(v\in S).
\]
To apply the bound obtained from fractional spanning-tree packing,
we check that these required values sum to an integer. Since
$\sum\limits_{v\in S}\eta(v)=\partial h'(s)$, we have
\[
\begin{aligned}
 \sum\limits_{v\in S}\beta(v)
 &=\frac{\sum\limits_{v\in S}d_G^-(v)}{2}-\partial h'(s)\\
 &\equiv\frac{d_{G/F}^-(s)}{2}-\partial h'(s)
 \equiv0\pmod1.
\end{aligned}
\]
Here $d_{G/F}^-(s)$ is the number of negative edges incident with
$s$ in the quotient. The first congruence holds because each negative
edge inside $F$ is counted twice in $\sum\limits_{v\in S}d_G^-(v)$,
whereas each negative edge leaving $S$ is counted once. The second
congruence is the required boundary condition for $h'$ at $s$.
Thus $\sum\limits_{v\in S}\beta(v)\in\Z$. By
Corollary~\ref{cor:prescribed-boundary} and the fractional spanning-tree
packing bound, there exists an assignment $h_F$ on $F$ satisfying
\[
 \partial h_F(v)\equiv\beta(v)\pmod1\quad(v\in S),
 \qquad |h_F(e)|\le\frac1{2\tauf(F)}<\frac1q\quad(e\in E(F)).
\]
Together with the retained values of $h'$ outside $F$, this gives
the required assignment on $G$, a contradiction.

For an arbitrary connected signed graph and arbitrary $q$, replace
each edge by two parallel copies with the same reference direction.
Assign signs to the two copies so that their product equals the sign
of the original edge.
Then the doubled graph has even degrees, each cut retains its sign,
and every cut has twice its original size. Applying the even-degree case
with parameter $2q$ gives an assignment with absolute value below
$1/(2q)$ on every copied edge and the boundary
congruences determined by the copied signs. Adding the values on the
two copies of each edge preserves the boundary and yields the
original congruences, with absolute value below $1/q$ on every original edge.
\end{proof}

\section{Circular flow consequences}\label{sec:flows}

We are now ready to prove Theorem~\ref{thm:flow}. We recall its statement.

\begin{flowtheoremrestated}
For every integer $q\ge3$, if every negative cut of a signed
graph $\widehat G$ has size at least $q$, then
$\Phi_c(\widehat G)<2q/(q-2)$.
\end{flowtheoremrestated}

\begin{proof}
We fix an arbitrary reference orientation of $G$.
We write $\nu:E(G)\to\{0,1\}$ for the additive encoding of the signature
$\sigma$, assigning $1$ to negative edges and $0$ to positive edges.
Theorem~\ref{thm:general-boundary}, applied to $\widehat G$, provides
an assignment $h$, which may be chosen rational as explained in
Remark~\ref{rem:quantitative-flow-gap}. We put $a=\max_e|h(e)|<1/q$ and define
\[
 r=\frac{2}{1-2a}<\frac{2q}{q-2}.
\]
We define $b(e)\in[0,r)$ by $b(e)\equiv r(h(e)+\nu(e)/2)\pmod r$.
Since $\partial\nu(v)\equiv d^-(v)\pmod2$ and
$\partial h(v)\equiv d^-(v)/2\pmod1$,
we have $\partial b(v)\equiv0\pmod r$ at every vertex $v$. As $|rh(e)|\le r/2-1$, the values $b(e)$ satisfy the required
edge bounds for both positive and negative edges.
Thus $b$ is a circular $r$-flow.
\end{proof}

\begin{remark}\label{rem:quantitative-flow-gap}
We record three observations.
\begin{enumerate}
\item The assignment $h$ in Theorem~\ref{thm:general-boundary} may be
chosen rational. Consequently, $a$ and $r$ in the flow construction
may also be chosen rational.
\item If the graph considered in Theorem~\ref{thm:flow} has $m\ge1$ edges,
the upper bound can be improved to
\[
 \Phi_c(\widehat G)\le\frac{2qm}{m(q-2)+1}
 =\frac{2q}{q-2}-\frac{2q}{(q-2)(m(q-2)+1)}.
\]
\item For every connected loopless signed graph $\widehat G$ with
$|V(G)|\ge2$ and $\tauf(G)>1$, we have
\[
 \Phi_c(\widehat G)\le\frac{2\tauf(G)}{\tauf(G)-1}.
\]
This follows from Corollary~\ref{cor:prescribed-boundary} and the
fractional spanning-tree packing bound. Equality can occur.
\end{enumerate}

\end{remark}

\section{Proof of the coloring theorem}\label{sec:duality}

\begin{proof}[Proof of Theorem~\ref{thm:main}]
We fix a plane embedding of $\widehat G$. Its signed dual
$\widehat G^*$ is the usual plane dual, with each dual edge receiving
 the sign of its corresponding primal edge.
We orient the primal edges arbitrarily and each dual edge from the
left face to the right face of its oriented primal edge.
Under plane duality,
switching corresponds to inversing on dual cycles, and
$g_-(\widehat G)\ge\ell$ is equivalent to every negative cut of
$\widehat G^*$ having size at least $\ell$.

By Theorem~\ref{thm:flow}, $\widehat G^*$ admits a circular $r$-flow
for some $r<2\ell/(\ell-2)$. The dual mapping is then a circular $r$-tension on $\widehat G$. Therefore, $\chi_c(\widehat G)\le r<2\ell/(\ell-2)$.
\end{proof}

\section{Tightness}\label{sec:tightness}

\subsection{Circular flows and modulo orientations}

Observing that the complete graph $K_{4k}$ has edge connectivity
$4k-1$ and admits no circular $(2+1/k)$-flow, Jaeger conjectured
that edge connectivity $4k$ suffices~\cite{Jaeger84}.
Zhang proposed that odd edge connectivity can replace edge connectivity
in this conjecture~\cite[Conjecture~5.2]{Zhang02}.
A major step toward Jaeger's conjecture was taken in 2013 by
Lov\'asz, Thomassen, Wu, and Zhang~\cite{LTWZ13}, who proved that every
graph of odd edge connectivity at least $6k+1$ admits a circular
$(2+1/k)$-flow.

In 2018, Han, Li, Wu, and Zhang~\cite{HLWZ18} constructed
$4k$-edge-connected counterexamples for every $k\ge3$, and
$(4k+1)$-edge-connected counterexamples for every $k\ge5$.
Thus, for every $k\ge5$, Theorem~\ref{thm:flow} gives the best
possible edge-connectivity threshold, namely $4k+2$, for
guaranteeing a circular $(2+1/k)$-flow.

In terms of modulo orientations, Theorem~\ref{thm:flow} and the
counterexamples above determine the optimal edge-connectivity
threshold for a modulo $(2k+1)$-orientation for every $k\ge5$.
The remaining cases are whether edge connectivity $4$ or $5$
guarantees a modulo $3$-orientation, $8$ or $9$ guarantees a
modulo $5$-orientation, $13$ guarantees a modulo $7$-orientation,
and $17$ guarantees a modulo $9$-orientation.
These optimal connectivity thresholds do not, however, determine
the best possible upper bound on the circular flow number for
each edge connectivity. Determining that bound remains open.

\subsection{Coloring and homomorphisms}
\paragraph{Simple planar signed graphs.}
Applying Remark~\ref{rem:quantitative-flow-gap}, we obtain, for every
simple planar signed graph $\widehat G$ on $n\ge3$ vertices,
\[
 \chi_c(\widehat G)\le6-\frac{6}{n-1}.
\]
Indeed, we may extend $\widehat G$ to a signed plane triangulation
on the same vertex set. Its dual $H$ is cubic and $3$-edge-connected,
with $2n-4$ vertices, and the fractional spanning-tree packing
formula gives $\tauf(H)=(3n-6)/(2n-5)$.

Wang and Zhou~\cite{WangZhou26} recently constructed a simple
planar signed triangulation on $85$ vertices with circular
chromatic number exactly $5$, improving the previous lower bound
of $14/3$. Thus the supremum of the circular chromatic numbers
of simple planar signed graphs lies between $5$ and $6$;
its exact value remains open.

\paragraph{Negative girth $2k$.}
Suppose that $\widehat G$ is a connected signed plane graph on
$n$ vertices, with negative girth at least $2k$, where $k\ge2$,
and every face is a cycle of length $2k$. Euler's formula gives
$|E(G)|=k(n-2)/(k-1)$. The edge-count bound in
Remark~\ref{rem:quantitative-flow-gap} therefore yields
\[
 \chi_c(\widehat G)\le \frac{2k}{k-1}
 -\frac{2k}{(k-1)(2k(n-2)+1)}.
\]
If, in addition, the ordinary girth is at least $2k$, the dual $H$
is $2k$-regular and $2k$-edge-connected. With
$f=|V(H)|=(n-2)/(k-1)$, the fractional spanning-tree packing
formula gives $\tauf(H)=kf/(f-1)$. The packing bound improves
 the estimate to
\[
 \chi_c(\widehat G)\le \frac{2k}{k-1}
 -\frac{2k}{(k-1)(n-1)}.
\]

Both improvements depend on the number of vertices and tend to
zero as $n$ tends to infinity. For even $k$, the bound $2k/(k-1)$ cannot be improved, as follows
from the construction in~\cite{NWZ21} by signed subdivision.

For $k=2$, the facial condition makes $G$ bipartite. For simple
bipartite planar signed graphs on $n\ge4$ vertices,
Kardo\v{s}, Narboni, Naserasr, and Wang~\cite{KNNW23} proved the
stronger bound
\[
 \chi_c(\widehat G)\le4-\frac{4}{\lfloor(n+2)/2\rfloor},
\]
which is sharp for every $n\ge4$.

\paragraph{Negative (odd) girth $2k+1$.}
For a signed planar graph $\widehat G$ with negative girth at least
$2k+1$, Theorem~\ref{thm:main} gives
\[
 \chi_c(\widehat G)<\frac{4k+2}{2k-1}=2+\frac{4}{2k-1}.
\]
For a graph with $m\ge1$ edges, Remark~\ref{rem:quantitative-flow-gap}
refines this to
\[
 \chi_c(\widehat G)\le
 \frac{2(2k+1)m}{(2k-1)m+1}.
\]

Planar graphs of odd girth $4k-1$ that do not admit a homomorphism
to $C_{2k+1}$ are obtained by taking the plane duals of the graphs
$W(2k-1)$ in~\cite[Section~5]{LTWZ18}.
These are subdivided wheels, a construction communicated by
DeVos~\cite{DeVos00}. Earlier examples are recorded in~\cite{GGH01}.

Corollary~\ref{cor:ordinary-target} shows that every planar graph
of odd girth at least $4k+3$ admits a homomorphism to $C_{2k+1}$.
The examples above show that odd girth $4k-1$ does not suffice.
Whether odd girth at least $4k+1$ suffices remains open for
$k\ge2$~\cite{Zhang02}. The case $k=1$ follows from
Gr\"otzsch's theorem, but is not implied by our general bound.

In recent joint work with Chen, MacGillivray and Nasserasr, we construct
 a $K_5$-minor-free graph of odd girth $6k-3$ (for every $k>1$)
that admits no homomorphism to $C_{2k+1}$.
This highlights the importance of planarity in this context.

\section*{Acknowledgment}
This work has received support under the program ``Investissement d'Avenir''
launched by the French Government and implemented by ANR, with the reference
``ANR-18-IdEx-0001'' as part of its program ``Emergence''.\\

The development and preparation of this work benefited from the assistance of ChatGPT through a paid subscription.
\begingroup
\small\raggedright
\bibliographystyle{unsrt}
\bibliography{Jaeger-Zhang-references}
\endgroup
\end{document}